\documentclass[12pt]{article}
\usepackage[mathscr]{eucal}
\usepackage{amssymb,amsmath}
\font\cmc=cmcsc10  scaled \magstep2

\newcommand\vk{\vskip}

\newcommand\al{\alpha}

\newcommand\la{\lambda}
\newcommand{\nhd}{neighborhood}
\newcommand\ga{\gamma}

\newcommand{\ve}{\varepsilon}
\newcommand{\va}{\varphi}

\newcommand{\Hol}{\text{Hol}}
\newcommand{\psh}{plurisubharmonic}

\newcommand\Om{\Omega}

\newcommand\de{\delta}
\newcommand\De{\Delta}

\newcommand\no{\noindent}

\newcommand\ovl{\overline}

\newtheorem{Proof.}{\it Proof.}

\begin{document}
\vbox to .5truecm{}
\large
\begin{center}
\cmc Hyperbolicity and Vitali properties of unbounded domains in Banach spaces
\end{center}
\centerline {Nguyen Quang Dieu, Nguyen Van Khiem and Le Thanh Hung}
\vk.5cm
{\normalsize
\no{\bf Abstract.} Let $\Om$ be a unbounded domain in a Banach space. In this work, we wish to impose {\it local conditions} on boundary point of $\Om$ (including the point at infinity) that guarantee complete hyperbolic of $\Om.$
We also search for local boundary conditions so that Vitali properties hold true for $\Om.$
These properties might be considered as analogues of the taut property in the finite dimensional case.
\vk.3cm
\no {\bf Keywords} Lempert functions, Kobayashi pseudo-distance, hyperbolic domains , complete hyperbolic domains

\no {\bf Mathematics Subject Classification 2010:} 32A10, 32A07, 32C25, 32Q45, 32U05 
}
\vk.5cm\no{\textbf {1. Introduction}}
\vk0,3cm
\noindent
Given a domain (open connected subset) $\Om$ in a Banach space, we say that $\Om$ is (Kobayashi) {\it hyperbolic} if 
the Kobayashi pseudo-distance is a distance and defines the topology of the domain.
Moreover, if every Cauchy sequence with respect to the Kobayashi pseudo-distance in $\Om$ is convergent 
then $\Om$ is called {\it complete hyperbolic}. 
On the other hand, 
we introduce in \cite{KDK} (see also \cite{HQVH}) some properties which are intermediate between hyperbolicity and complete hyperbolicity. Loosely speaking, $\Om$ is said to have the Vitali property if 
every sequence of holomorphic mappings from a domain $A$ in some Banach space into $\Om$ that converges 
pointwise only on a sufficiently large subset of $A$ must converge locally uniformly.
For the ease of the exposition, $\Om$ is said to have the weak (resp. strong) Vitali property (WVP) (resp. (SVP))
depending on the nature of the locus where the pointwise convergence holds.
We must say that the above Vitali properties are inspired from a classical theorem due to Vitali which asserts that a locally uniformly bounded sequence of scalar holomorphic functions from a domain $D$ in $\mathbb C$ is locally uniformly convergent if it is pointwise convergent on a set having an accumulation point in $D$.
Unfortunately, we do not know if SVP and WVP are really different properties.
We are, however, able to show in \cite{KDK} the following implications: 
$$\text{complete hyperbolic}\  \Rightarrow \text{SVP} \Rightarrow \text{WVP} \Rightarrow \text{hyperbolic}.$$
Moreover, it is shown in Theorem 3.14 of \cite{KDK} that in the case of domains in $\mathbb C^n$, the two properties SVP and WVP agree and coincides with the notion of {\it tautness}. Recall that a domain $\Om \subset \mathbb C^n$ is called taut if each sequence 
$\{f_j\} \subset \Hol(\De,\Om)$ contains a subsequence which is either convergent or compactly divergent.
It was observed in \cite{HQVH} (see also \cite{KDK}) that the above notion of tautness does not admit
a natural generalization to the infinite dimensional case. Indeed, consider the sequence 
$$f_j (\la):=(0,\cdots,0, \frac{\la}2,0,\cdots,0, \cdots), \ \la \in \De,$$
where $\De$ is the unit disk in $\mathbb C.$
Then $\{f_j\}$ is a sequence of holomorphic mappings from $\De$ into the unit ball of $c_0$
that contains no subsequence which is either convergent or compactly divergent. 
It means that the unit ball of $c_0$ does not have the expected taut property. According 
to \cite{KDK} (see also \cite{HQVH}) the Vitali properties seem to be appropriate analogues of the taut property in infinite dimensional spaces. 

A common theme in complex analysis is to decide whether a domain has some property if locally near each boundary point it has. In the finite dimensional case, we know that a {\it bounded} domain $\Om$ is taut or compete hyperbolic if 
and only if $\Om$ is locally taut or locally complete hyperbolic, respectively.
The aim of this note is to extend these results to the case of {\it unbounded} domains in Banach spaces. 
Here is a brief outline of our work. In the next section, we review basic elements of complex analysis in Banach spaces pertaining to our work. We recall, among other things, the Lempert functions and the Kobayashi pseudo-distance. 
Their construction are completely analogous as in the finite dimensional case. Nevertheless, these concepts, as in \cite{NP}, allow us to introduce the key notions of $k'-$ and $t-$boundary points.
Our main results are explained in Section 3. Theorem 3.1 characterizes hyperbolicity of unbounded domains in Banach spaces in terms of the behavior at infinity of the Lempert function. In the finite dimensional case, this result is exactly Proposition 3.1 in \cite{NP}. We should say that the proof given there does not extend directly to our case since it requires local compactness of the ambient space. 
The same remark applies to our next main result, Theorem 3.5, which in the same spirit as Proposition 3.6 of \cite{NP}
relates complete hyperbolicity of a domain $\Om$ and the $k'-$property of its boundary points.
Moreover, even in the finite dimensional case, our proof is perhaps simpler and more direct than the original proof
in \cite{NP}. Namely, we rely essentially on a comparison principle for Kobayashi pseudo-distances on domains in Banach spaces (see Proposition 3.6). We continue our investigation by presenting in Theorem 3.8 a result which says that SVP is indeed a {\it local} property for domains having $\infty$ as a $t-$point. This result is, again, an infinite dimensional version of Proposition 3.9 in \cite{NP} in which the taut property is now replaced by SVP.
The final result Proposition 3.11 allows us to cook up more example of domains having SVP by using holomorphic transformations between domains in Banach spaces. This type of result, in the finite dimensional case,  is only discussed in Proposition 5.1.8 of \cite{Ko}.
We also provide explicit applications of the above mentioned results in the remarks following them.
\vskip1cm
\noindent
{\bf Acknowledgments.} We are greatly indebted to Professor Nguyen Van Khue for sharing with us his deep insight and  useful suggestions during the preparation of this manuscript. This work would not have been written without his generous help. The first named author is supported by the grant 101.02-2016.07 from the NAFOSTED program. 
\vk0,5cm 
\no{\textbf {2. Preliminaries}}
\vk.5cm
\noindent
We first introduce some standard notation.
Given open connected subsets $A, \Om$ of Banach spaces, we write $\Hol(A,\Om)$ for the space of holomorphic mappings 
from $A$ into $\Om.$ This space is equipped with the topology of local uniform convergence.
We also let $\De (0,r):=\{z \in \mathbb C: \vert z\vert<r\}$ and $\De:=\De(0,1)$ for short.
More generally, by $\mathbb B(a,r)$ we mean the open ball with center $a$ and radius $r>0$ in some Banach space.

Next, we move to the construction of the Lempert function, the first fundamental concept in our work.
To this end, we require the following easy fact whose proof is left to the reader.

\noindent
{\bf Lemma 2.1.} {\it Let $\al, \beta \in \De, \al \ne \beta$ and 
$\al', \beta' \in E$. Then there exists $\la_1, \la_2 \in E$ such that for $h(t)=\la_1 t+\la_2 (t \in \mathbb C)$
we have
$$h(\al)=\al', h(\beta)=\beta', \Vert h\Vert_{\De} \le \Vert \al'\Vert +2\Vert \frac{\al'-\beta'}{\al-\beta}\Vert.$$}
\noindent
The following result clarifies the arguments given in Remark 3.1.1 of \cite{JP}.
\vskip0,2cm
\noindent
{\bf Lemma 2.2.} {\it Let $\Om$ be a domain in a Banach space $E$ and $z, w \in \Om.$ Then there exists
$f \in \Hol(\De, \Om)$ such that $z, w \in f(\De).$}
\vskip0,2cm
\noindent
{\it Proof.} First, since $D$ is a connected open set, we can choose a continuous map $\gamma: [0, 1] \to \Om$
such that $\gamma([0,1]) \Subset \Om, \gamma(0)=z, \gamma (1)=w.$
Choose $\ve>0$ so small such that
$$\gamma ([0,1])+\bar {\mathbb B}(0,4\ve) \subset \Om.$$
Next, we consider a sequence of Bernstein polynomials
$$P_n (t):= \sum_{k=0}^n \gamma(\frac{k}{n}) C^k_n t^k(1-t)^{n-k}, t \in [0,1].$$ 
As in the case where $E=\mathbb R$, we can show that $P_n$ converges uniformly to $\gamma$ on $[0,1].$
In particular,
there exists a polynomial $P: \mathbb C \to E$ such that $\Vert P-\ga\Vert_{[0,1]}<\ve.$
By Lemma 2.1, there exists a polynomial $h$ of degree $1$ in $\mathbb C$ with values in $E$ that satisfies

\noindent
(a) $h(0)=z-P(0), h(1)=w-P(1).$

\noindent
(b) $\Vert h \Vert_\De \le \Vert z-P(0) \Vert+\Vert (z-P(0))-(w-P(1))\Vert<3\ve.$

Set 
$$g(t):=P(t)+h(t), \ t \in \De.$$
It follows from (a) and (b) that $g(0)=z, g(1)=w$ and
$$\Vert g- \ga\Vert_{[0,1]}<4\ve.$$
Therefore $g([0,1]) \Subset \Om,$ by the choice of $\ve.$
Finally, we let $U$ be bounded simply connected domain in $\mathbb C$ such that $[0,1] \subset U$ and
$g(U) \subset \Om.$
Let $\va: \De \to U$ be a biholomorphic mapping and $f:= g \circ \va \in \Hol(\De,E)$. Then we have 
$$z, w \in f(\De)=g(U) \Subset \Om.$$ 
The proof is thereby completed.
\vk0,2cm
\noindent 
In virtue of the above result, it makes sense to define the Lempert function for an arbitrary domain $\Om$ as follows
$$\begin{aligned} l_\Om (z, w)&:=\inf \{\vert \la\vert: \exists f \in \Hol (\De, \Om), f(\la)=z, f(0)=w\}\\
&=\inf \{\vert \la\vert: \exists f \in \Hol (\De,\Om), f(0)=z, f(\la)=w\}.
\end{aligned}$$
One advantage of the Lempert functions is the following decreasing property: If $F: \Om \to \Om'$ is a holomorphic mapping between domains in Banach spaces then
$$l_{\Om'} (F(a), F(b)) \le l_\Om (a,b), \ \forall a,b \in \Om.$$
This fact follows immediately from the definition of Lempert functions. A less obvious property is provided by the next result whose proof is the same as the one given in the finite dimensional case (see Proposition 3.1.13 in \cite{JP}). 
\vk0,2cm
\noindent
{\bf Lemma 2.3.} {\it The function $l_{\Om}(z,w)$ is upper semicontinuous on $\Om \times \Om.$}
\vk0,3cm
\noindent
Following \cite{NP}, we will use Lempert functions to analyze boundary points (possibly at infinity) of a domain $\Om$ in a Banach space $E$. More precisely, 
we say that $a \in \partial \Om \cup \{\infty\},$ is a (global) $t-$point of $\Om$ if
$$\lim_{z \to a, w \to b} l_\Om (z,w)=1, \ \forall b \in \Om.$$
\noindent
We now turn to the construction of the Kobayashi pseudo-distance on a domain $\Om$ in a Banach space.
For $\la \in \De$ we let 
$$p(\la):=\frac1{2}\log \frac{1+\vert \la\vert}{1-\vert \la\vert},$$ 
be the Poincare distance from $0$ to $\la.$ As in the finite dimensional case (see \cite{JP}, p.73), the Lempert function can now be used to define the Kobayashi pseudo-distance on $\Om$ as follows. 
For $z,w \in \Om$ we let
$$k_\Om (z,w):=\inf\{\sum_{i=1}^n p(l_\Om (z_{i-1},z_i)): z_0=z, z_1,\cdots, z_{n-1} \in \Om, z_n=w\}.$$
It is easy to see that $k_\Om$ defined as above coincides with the original construction (see \cite{Ko}, p.50 or \cite{FV} p.81) using holomorphic chains. 

By the same proof as in the finite dimensional case (see Proposition 3.1.7 in \cite{Ko}), we can show that 
$k_\Om$ is decreasing under holomorphic maps i.e, if $f: \Om \to \Om'$ is a holomorphic mapping between domains in Banach spaces then
$$k_{\Om'} (f(z),f(w)) \le k_\Om (z,w),\ \forall z, w \in \Om.$$
Then we say that $\Om$ is {\it hyperbolic} if $k_\Om$ is a distance and defines the topology of $\Om$.
Notice that, in contrast to the case finite dimensional case, $k_\Om$ may be a distance without defining the topology of $\Om$ (see \cite{FV} p. 93).
Furthermore, $\Om$ is said to be {\it complete hyperbolic} if every $k_\Om-$ Cauchy sequence in $\Om$ is convergent.
By Proposition 6.9 in \cite{D} (see also Proposition 3.6 in \cite{HQVH}) we know that every bounded convex domain $\Om$ in a Banach space is complete hyperbolic.
Hence, all connected open subsets of $\Om$ are hyperbolic. In particular, each bounded connected open subset of a Banach space is hyperbolic.

The next ingredient needed in our work is the concept of plurisubharmonic functions on Banach spaces.
More precisely, we say that $\va: \Om \to [-\infty, \infty),$ where $\Om$ is a domain in a Banach space, is plurisubharmonic if $u$ is upper semicontinuous on $\Om$ and the restriction of $u$ on the intersection of $\Om$ with each complex line in $E$ is subharmonic. Notice that we allow the function $u \equiv -\infty$ to be plurisubharmonic. 

Using the Kobayashi pseudo-distance and plurisubharmonic functions, we obtain below a sufficient condition for a boundary point to be a $t-$point. This slightly generalizes Proposition 3.4 in \cite{NP}.
\vk0,2cm
\noindent
{\bf Proposition 2.4.} {\it Let $\Om$ be a  domain in a Banach space $E$ and $a \in \partial \Om \cup\{\infty\}$. Assume that there exists $u \in PSH(D), u<0$ satisfying the following conditions:

\noindent
(a) $u$ is a barrier at $a$ i.e., $\lim_{z \to a} u(z)=0.$ 

\noindent
(b) If $\{x_n\} \subset \Om$ and $u(x_n) \to 0$ then $\{x_n\}$ does not converge in $k_\Om.$

Then $a$ is a $t-$point of $\Om.$}

\noindent
Note that, in the case where $\Om$ is hyperbolic, the condition (b) is trivially satisfied.
\vk0,2cm
\noindent
{\it Proof.} Assume that $a$ is not a $t-$point. Then there exist $\{f_j\} \subset \Hol(\De,\Om)$ 
and $\{\la_j\}\subset \De$ with $f_j (0) \to b \in \Om, f_j(\la_j) \to a$ and $\la_j \to \la_0 \in \De$.
Set $$v(\la):= \sup_{j \ge 1} (u \circ f_j)(\la), \ \forall \la \in \De.$$
Then $v^*$, the upper regularization of $v$, is $\le 0$ and subharmonic on $\De$. Moreover, by the assumption (a) we obtain
$$v(\la_j ) \ge (u \circ f_j)(\la_j) \to 0 \ \text{as}\ j \to \infty.$$
It follows that $v^*(\la_0)=0.$
The maximum principle now implies that 
$v^* \equiv 0$ on $\De$. Since $v=v^*$ almost everywhere on $\De$, we may choose $\{\beta_n\} \subset \De$ with 
$\beta_n \to 0$ and $v(\beta_n)=0$. So for each $n \ge 1$ we can choose $j_n$ such that 
$$(u \circ f_{j_n})(\beta_n)>-\frac1{n}.$$
Set $x_n:=f_{j_n}(\beta_n)$. We have $u(x_n) \to 0$ as $n \to \infty$. By the triangle inequality and the decreasing property of the Kobayashi pseudo-distance we get
$$k_\Om (x_n,b) \le k_\Om (x_n,f_{j_n}(0))+k_\Om (f_{j_n}(0),b) \le k_\De (\beta_n, 0)+k_\Om (f_{j_n}(0),b), \ \forall n \ge 1.$$ 
By letting $n \to \infty$ we have $k_\Om (x_n, b) \to 0$. This is a contradiction to (b). The desired conclusion now follows.
\vskip0,2cm
\noindent
The following notions are natural generalizations of the corresponding ones in the case of finite dimensional case
(see \cite{NP}, p.611).
\vk0,2cm
\noindent
{\bf Definition 2.5.} {\it Let $\Om$ be a domain in a Banach space $E$ and 
$a \in \partial \Om \cup \{\infty\}.$ We say that $a$ is a $k'-$point of $\Om$ if there is no
$k_\Om-$Cauchy sequence of $\Om$ that converges to $a$. The point $a$ is called a local $k'-$ point of $\Om$ if it admits  a \nhd\ $U$ such that $a$ is a $k'-$point for each connected component of $U \cap \Om.$} 
\vk0,2cm
\noindent
The above notion will play a central role in our investigation of complete hyperbolic domains in Banach spaces.
The relationship between $t-$points and $k'-$points is, however, not clear to us.
We rely on the following result on checking whether a boundary point of a domain in Banach space is a $k'-$point (resp. local $k'-$point).
\vk0,2cm
\noindent
{\bf Proposition 2.6.} {\it Let $\Om$ and $\Om'$ be domains in Banach spaces and 
$\va: \Om \to\Om'$ be a holomorphic mapping. Let $a \in \partial \Om \cup \{\infty\}$ be a boundary point having the following property: For every sequence $\{a_n\} \subset \Om, a_n \to a$, there exists a subsequence $\{a_{n_k}\}$ such that 
$\va(a_{n_k})$ converges to $a' \in \partial \Om' \cup \{\infty\}$ which is a $k'-$ point (resp. local $k'-$point) of $\Om'.$

Then $a$ is a $k'-$point (resp. local $k'-$point) of $\Om.$}
\vk0,2cm
\noindent
The proof of this proposition is a straightforward application of the  decreasing property of the Kobayashi pseudo-distance. The details are therefore omitted.

Finally, we recall the notions of Vitali properties inverstigated in \cite{KDK}. Before giving the precise definitions, 
we introduce the following notation: Given a subset $S$ of a domain $A$ in a Banach space, we let

\no
$S^u:=\Big\{z\in A \cap\ \ovl S: \forall$ connected neighborhood $U$ of $z$ and every 
$f \in \Hol (U, \mathbb C), f\big|_{ U\cap S}=0 \Rightarrow f\big|_ U=0  \Big\}$.
\vk0,2cm
\noindent
{\bf Definition 2.7.} {\it Let $\Om$ be a domain in a Banach space. Then we say that $\Om$ has the weak Vitali property 
(WVP for short) if for every sequence of holomorphic mappings $\{f_j\}$ from a connected open set $A$ in a Banach space
into $\Om$ is convergent in $\Hol (A,\Om)$ provided that $Z_{\{f_j\}} \cap Z^u_{\{f_j\}} \ne \emptyset$.
Here $Z_{\{f_j\}}$ denote the collection of points $x \in A$ such that $f_j (x)$ is convergent.
Moreover, $\Om$ is said to have the strong Vitali property (SVP for short) if the above sequence $\{f_j\}$ is convergent
in $\Hol (A,\Om)$ as long as $Z^u_{\{f_j\}} \ne \emptyset$.
}
\vk0,2cm
\noindent
In \cite{KDK}, these properties are formulated in the broader context of Banach analytic manifolds.
Nevertheless, we can show, in the case of domains in Banach spaces, Definition 2.7 agrees with Definition 2.1 in \cite{KDK}.
At the end of this paper, we will prove that the two Vitali properties again coincide in the category of 
unbounded domains having $\infty$ as a $k'-$point.
\vk0,5cm
\no{\textbf {3. Main Results}}
\vk.3cm
\noindent
The first result generalizes Proposition 3.1 in \cite{NP} to the case of unbounded domains in Banach spaces.

\no {\bf Theorem 3.1.} {\it Let $\Omega$ be a unbounded domain in a complex Banach space. Then the following statements are equivalent:

\noindent
(a) $\Omega$ is hyperbolic. 

\noindent
(b) $\liminf\limits_{z\to\infty,\ w\to b }\ell_{\Omega}(z, w)>0\quad\forall\ b\in \Omega.$}
\vk0,2cm
\no Since Montel's theorem is not valid for Banach-valued holomorphic functions, the proof given in Proposition 3.1 in \cite{NP} does not directly apply in our case. Instead, we employ ideas from Theorem 3.1 in \cite{KDK} that relates hyperbolicity and Vitali property of domains in Banach spaces.
We recall the following auxiliary results which are taken from \cite{KDK}. The first lemma is an analogue of an earlier result due to Kiernan proved in the finite dimensional case (see Lemma 5.1.4 in \cite{Ko}).
\vk.3cm
\no{\bf Lemma 3.2} (\cite{KDK}, Lemma 3.2). {\it Let $Y$ be a connected open subset of a complex Banach space and $x \in Y.$ Let $U, V, W$ be open subsets of $Y$ such that
$x \in V \Subset U, \overline{U} \cap \overline{W}=\emptyset$ and $U$ is hyperbolic.
Assume that there exists $\de \in (0, 1)$ such that for every  $f \in \Hol (\De, Y)$ with $f(0) \in V$ we have 
$f(\De(0,\de)) \subset U$. Then $k_Y (x, W)>0.$}
\vk.3cm 
\no
We also need the following variant of Vitali's convergence theorem for holomorphic vector-valued functions. The proof is  a slight modification of the proof of Theorem 2.1 in \cite{AN}.
\vk0,2cm
\noindent
{\bf Lemma 3.3} (\cite{KDK}, Lemma 3.5). {\it Let $E, F$ be Banach spaces and $\Omega$ be an open subset of $E.$ Let $\{f_j\}$ be a sequence in $\text{Hol}(\Om, F)$ that satisfies the following conditions:

\noindent
(a) $\{f_j\}$ is locally uniformly bounded on $\Om.$

\noindent
(b) $Z_{\{f_j\}}$ is a set of uniqueness for $\Hol(\Om, \mathbb C)$ i.e., every holomorphic function $g: \Om \to \mathbb C$ that vanishes on $Z_{\{f_j\}}$ must be identically $0.$

\no Then $\{f_j\}$ converges in $\text{Hol}(\Om,F)$.}
\vk.3cm 
\no{\it Proof of Theorem 3.1.} $(a) \Rightarrow (b).$ Assume that there exists $b\in \Omega$ such that
$$\liminf\limits_{z \to\infty,\ w\to b }\ell_{\Omega}(z, w)=0.$$
Then there exists a sequence $\{f_j\}\subset \text{Hol}(\Delta,\Omega)$ such that 
$f_j(0) \to \infty, f_j(x_j)\to b$ with $x_j\to 0$ and $x_j\in\Delta$. Since $\Omega$ is hyperbolic, by the decreasing distance property of Kobayashi pseudo-distance we obtain
$$k_\Om (f_j (0), f_j (x_j)) \le k_\De (0, x_j) \to 0\;\;\text{as}\;\; j\to \infty.$$
It follows that $f_j(0)\to b$ as $j\to \infty$. This is impossible.
\vk.3cm 
\noindent
$(b) \Rightarrow (a).$ First we note that if $x_n \to x$ in the original topology of $\Om$ then 
$k_\Om (x_n, x) \to 0$ as $n \to \infty.$
Conversely, assume that there exists a sequence $\{x_n\} \subset \Om$ such that
$k_{\Om}(x_n, x)\to 0$ as $n\to\infty$ but $x_n\nrightarrow x$. Take a bounded open neighborhood $U$ of $x$ and an open neighborhood $W$ of the sequence $\{x_n\}$ such that $\overline U\cap \overline W=\emptyset$. Then $U$ is hyperbolic. 
Choose a sequence of open sets $V_n$ such that
$U\supset V_n \downarrow x$. Define a sequence$\{\delta_n\}_{n\geq 0}$ by $\delta_0=\dfrac12$, $\delta_1=\dfrac13$ and
$$0<\de_{n+1}<\min \Big \{\frac1{n},\ r_n:=\de_n \prod_{j=0}^{n-1}\frac{\de_j-\de_n}{1-\de_j\de_n}\Big \}, \ \forall n \ge 1.$$
It follows that $r_{n+1}<\de_{n+1}<r_n.$ In particular, $r_n \downarrow 0$. Using Lemma 3.2, we obtain a sequence 
$\{f_n\} \subset \Hol(\De,\Omega)$ and points $a_n \in \De(0, r_n)$ such that
$$f_n (0) \in V_n, f_n(a_n) \not\in U, \ \forall n \ge 1.$$
For each $n\geq 1$, define
$$\theta_n(\lambda):=\frac{a_n}{r_n}\la \prod_{j=0}^{n-1}\dfrac{\de_j-\lambda}{1-\de_j\lambda},\  \forall \lambda\in\Delta.$$
Then $\theta_n\in \text{Hol}(\Delta, \Delta)$. Moreover, we have
$$\theta_n (0)=\theta_n (\de_j)=0, \ \forall 0 \le j \le n-1;\; \theta_n (\de_n)=a_n.$$
Finally, we define for each $n \ge 1$
$$g_n:=f_n\circ \theta_n\in \text{Hol}(\Delta, \Omega).$$
Then we have
$$g_n (\de_j)=g_n(0)=f_n(0), \ \forall n \ge 1, \forall 0 \le j \le n-1.$$
This implies that
$$\lim\limits_{n\to\infty}g_n(0)=\lim\limits_{n\to\infty}g_n(\de_j)=\lim\limits_{n\to\infty}f_n(0)=x \quad\forall j\geq 0.$$
$$g_n (\de_n)=f_n (a_n)\not\in U.$$
By the hypothesis, we have
$$r:=\liminf\limits_{z \to\infty,\ w\to x}\ell_{\Omega}(z, w)>0.$$
We claim that $\{g_n\}$ is locally bounded on $r\Delta$. Indeed, if not there exists $\lambda_n\in r\Delta$, $\lambda_n\to \lambda_0\in r\Delta$ such that $z_n=g_n(\lambda_n)\to\infty$. This implies 
$$\liminf\limits_{n\to \infty}\ell_{\Omega}(z_n, g_n(0))\leq |\lambda_0|<r.$$
Since $g_n (0)\to x$ we obtain a contradiction to the choice of $r.$
Now we use Lemma 3.3 to get that $g_n\to g$ in $\Hol(r\Delta , \Omega)$. This is impossible, because
$$a_n\to 0\;\; \text{and} \;\; g_n(a_n)\not\in U\ni x\;\;\forall n\geq 1.$$ 
The proof is complete.
\vskip0,2cm
\noindent
{\bf Remark.} If $\infty$ is a $k'-$point of $\Om$ then $\Om$ is hyperbolic. Indeed, if this is false 
then there exist $a_n \to \infty, b_n \to b\in \Om$ such that $l_\Om (a_n, b_n)\to 0.$ The triangle inequality now yields that
$$k_\Om(a_n, b)\le k_\Om (a_n, b_n)+k_\Om (b_n,b) \le p(l_\Om (a_n,b_n))+k_\Om (b_n,b) \to 0 \ 
\text{as}\ n \to \infty.$$
Hence $a_n \to \infty$ is a $k_\Om-$Cauchy sequence which is absurd. 
\vskip0,3cm
\noindent
Let $E, F$ be Banach spaces,
$\Omega$ be a domain in $E$ and $h:\Omega\times F\longrightarrow \mathbb R$ be a non-negative upper semicontinuous function satisfying 
$$h(z, \lambda w)=|\lambda|h(z, w), \  \forall \lambda\in\Delta.$$
Define
$$\Omega_h:=\{(z, w)\in \Omega\times F: h(z,w)<1\}.$$
Then $\Om_h$ is a Hartogs domain over $\Om$ with (balanced) fibers in $F$. 
Using Theorem 3.1, we are able to characterize hyperbolicity of $\Omega_h$ in terms of $\Om$ and $h.$
\vk.3cm 
\no{\bf Corollary 3.4}. {\it $\Omega_h$ is hyperbolic if only if $\Omega$ is hyperbolic and
$$\inf\limits_{K\times \partial\mathbb B}h(z, w)>0\;\;\forall K\;\text{compact in}\; \Omega,$$
where $\mathbb B=\{w\in F: \|w\|<1\}$.}
\vk.3cm
\no{\it Proof}. We follow the lines of the proof of Proposition 4.2 in \cite{NP}.
First, assume that $\Omega_h$ is hyperbolic. Since $\{\Om\} \times \{0\} \subset \Om_h,$
we infer that $\Omega$ is hyperbolic as well. Now, suppose that there exists $K\Subset \Omega$ such  that 
$$\inf\limits_{K\times \partial\mathbb B}h(z, w)=0.$$ 
Then there exists $a_j\in \Omega_h$, $a_j=(a_j', a_j'')\in K\times \partial \mathbb B$ such that $h(a_j)\to 0$. We may assume that $a_j'\to a'$. For each $j \ge 1$, we define 
$$f_j(t)=(a_j',\ t(h(a_j))^{-1})\in \Hol(\Delta, \Omega_h).$$
It is easy to check that
$$f_j(0)=(a_j',0)\to a^*:=(a', 0)\in K\times \{0\}\subset \Omega_h,$$
and for every $\de \in (0,1)$ we have
$$f_j(t)=(a_j',\ t(h(a_j))^{-1})\to \infty\;\;\forall t\in\Delta\setminus \De (0, \de).$$
Thus, for $\de \in (0,1),$ by the decreasing property of the Lempert functions we obtain
$$l_\Om (f_j(\de), f_j(0)) \le l_\De (\de, 0)=\de.$$
This implies that
$$\liminf\limits_{z\to \infty, w\to a^*}\ell_{\Omega_h}(z,w)=0,$$ which contradicts Theorem 3.1.
\vk.3cm 
\no Conversely, assume that $\Omega_h$ is not hyperbolic. By Theorem 3.1, there exists $\lambda_j\in \Delta$ and $f_j\in\Hol(\Delta, \Omega_h)$ with $\lambda_j\to 0$, $f_j(\lambda_j)\to\infty$ and $f_j(0)\to a\in \Omega_h$. Write $a=(a', a'')$ with $a'\in\Omega$, $a''\in F$ and $f_j=(f_j', f_j'')$ then $f'_j(\lambda_j)\to a'$, by hyperbolicity of $\Omega$. 
Set
 $K:= \{f'_j(\lambda_j)\} \cup \{a'\}$. 
Then $K$ is a compact subset of $\Om.$ Moreover, 
$f_j''(\lambda_j)\to\infty$ and we have
$$\begin{aligned}1>h(f_j(\lambda_j))&=\|f''_j(\lambda_j)\|h\Big(f_j'(\lambda_j), \dfrac{f''_j(\lambda_j)}{\|f_j''(\lambda_j)\|}\Big)\\ 
&\geq \|f''_j(\lambda_j)\|\inf\limits_{K\times \partial\mathbb B}h(z, w)\to\infty
\end{aligned}$$
which is a contradiction.
\vk.5cm 
\no
Our next main result characterizes complete hyperbolicity of unbounded domains in Banach spaces.

\noindent
{\bf Theorem 3.5}. {\it Let $\Omega$ be a unbounded domain in a Banach space such that $\infty$ is a $k'$-point of $\Omega$. Then the following assertions are are equivalent:

\no (a) $\Omega$ is complete hyperbolic;

\no (b)any finite boundary point of $\Omega$ admits an open neighbourhood $U$ of $p$ such that each connected component of $U\cap \Omega$ is complete hyperbolic;

\no (c) any finite boundary point of $\Omega$ is a local $k'$- point;

\no (d) any finite boundary point of $\Omega$ is a  $k'$- point.}
\vk.3cm 
\noindent
Before going into the proof, a few remarks are now in order.
\vk0,2cm
\noindent
{\bf Remarks.} (a) The main thrust of the theorem is the implication $(c) \Rightarrow (a)$ that characterizes complete hyperbolicity of $\Om$ in terms of {\it local} $k'-$ finite boundary points.

\noindent
(b) The derivation $(d) \Rightarrow (a)$ is false if $\infty$ is only supposed to be a {\it local} 
$k'-$point. For a trivial counterexample, we may take $\Om=\mathbb C.$ 

\noindent
(c) In the finite dimensional case, the above theorem is partially contained in Proposition 3.6 of \cite{NP}.
The proof given in Proposition 3.6 of \cite{NP} does not, however, directly generalize to our case since it requires the local compactness of $\mathbb C^n.$
See for example the implication $(iv) \Rightarrow (i)$ in \cite{NP} which corresponds to $(d) \Rightarrow (a)$ in our case.
\vskip0,2cm
\noindent
For the proof of Theorem 3.5, we first introduce the following notation.
Let $\Om$ be a domain in a Banach space and $a \in \Om.$ For each $\de>0$, we denote by $U_\Om (a,\de)$ the 
Kobayashi "ball" $\{x \in \Om: k_\Om (a,x)<\de\}.$ 
The needed technical result is the following comparison principle.
\vk0,2cm
\no{\bf Lemma 3.6.} {\it Let $\Om, a$ be as above and $\de, \ve$ be positive numbers. Then the following statements hold true:

\noindent
(a) $U_\Om (a,\de)$ is connected.

\noindent
(b) There exists a constant $C>1$ such that for every $\ve'>0$, each pair $p, q \in U_\Om (a,\de)$ can be joined by a chain of holomorphic disks lying in 
$U_\Om (a,3\de+\ve)$ with length not exceeding $C(k_\Om (p,q)+\ve').$ 
In particular,
$$k_{U_\Om (a,3\rho+\ve)} (p,q) \le C k_\Om (p,q), \ \forall p, q \in U_\Om (a,\rho).$$}
\noindent
{\it Proof.} (a) Fix $z \in U_\Om (a,\de).$ Since $k_\Om (a,z)=\de'<\de$, we may find a chain of points $z_0=a, z_1,\cdots, z_{n-}, z_n=z$ lying in $\Om$,
holomorphic maps $f_1,\cdots, f_n \in \Hol(\De,\Om)$ and points $a_1, \cdots, a_n$ of $\De$
such that
$$f_i(0)=z_{i-1}, f_i(a_i)=z_i, \ 1 \le i \le n$$
and
$$p(0,a_1)+\cdots +p(0,a_n)\le \de'.$$
Fix $1 \le i \le n$ and $\xi \in \De$ with $\vert \xi\vert \le \vert a_i\vert$. Then by the triangle inequality
and the decreasing property for the Kobayashi pseudo-distance we obtain
$$\begin{aligned}
k_\Om (a,f_i (\xi))&\le k_\Om (a,z_1)+\cdots+k_\Om(z_{i-1}, f_i(\xi))\\
&\le p(0,a_1)+\cdots+p(0,a_{i-1})+p(0,\xi) \le \de'<\de.
\end{aligned}$$
Thus $f_i (\xi)\in U_\Om(a,\de).$ This implies that $a$ and $z$ can be joined by a continuous curve sitting in side 
$U_\Om (a,\de).$ Hence $U_\Om (a,\de)$ is connected.

\noindent
(b) The proof follows exactly the same lines as in Proposition 3.1.19 in \cite{Ko}. 
The key idea is to express the Kobayashi pseudo-distance as the infimum of the lengths of holomorphic chains and then apply the decreasing property of the Kobayashi pseudo-distance. We do not repeat the details here.
\vk0,3cm
\no{\it Proof of Theorem 3.5.} $(a)\Rightarrow (b)\Rightarrow (c)$ are obvious.

\no
$(c) \Rightarrow (d).$ Assume for the sake of seeking a contradiction that there exists a $k_{\Omega}$-Cauchy sequence 
$\{z_n\} \subset \Omega$ such that $z_n \to a \in \partial \Om.$ For $n \ge 1$ we set
$$\rho_n:=\sup_{m \ge n, l \ge n} k_\Om (z_m, z_l).$$
Then $\rho_n \downarrow 0.$ We also denote for each $r>0$ the open set
$\Om_r:= \mathbb B(a, r)\cap \Om.$
We now claim that there exists $N \ge 1$ such that
$$U_\Om (z_N,4\rho_N)\subset \Om_N.$$
Indeed, if the claim is false then for each $n \ge 1$, there exists $w_n \in \Om$ such that 
$$\Vert w_n-a \Vert \ge n, k_\Om (z_n, w_n) <4\rho_n.$$
By the triangle inequality, we see that $\{w_n\}$ is a $k_{\Omega}$-Cauchy sequence. Since $w_n \to \infty$, we obtain a contradiction to the assumption that $\infty$ is a $k'-$point of $\Om.$ The claim now follows.
By Lemma 3.6 (a), there exists a connected component $\Om'_N$ of $\Om_N$ 
that includes $U_\Om (z_N, 4\rho_N).$
Moreover,
we can find a constant $C>1$ such that for $m \ge N, l \ge N$ we have
$$k_{\Om'_N} (z_m,z_l)\le k_{U_\Om (z_N, 4\rho_N)}(z_m, z_l) \le C k_\Om (z_m, z_l).$$
This implies that $\{z_n\}_{n \ge N} \subset \Om'_N$ is a $k_{\Om'_N}-$Cauchy sequence. 
Now we choose $r \in (0, N)$ such that $a$ is a $k'-$point for each connected component of $\Om_r.$  
For $n \ge 1$, we set
$$\eta_n:= \sup_{m \ge n, l \ge n} k_{\Om'_N} (z_m, z_l).$$
Then $\eta_n \downarrow 0.$ We now prove that there exists $N' \ge N$ such that
$$U_{\Om'_N} (z_{N'},4\eta_{N'}) \subset \Om_r.$$
Assume otherwise, then for each $n \ge N$, 
there exists $w_n \in \Om'_N \setminus \mathbb B(a,r)$ such that
$$k_{\mathbb B(a,N)}(z_n, w_n) \le k_{\Om'_N} (z_n, w_n) \le 4\eta_n \to 0 \ \text{as}\ n \to \infty.$$
This yields a contradiction to hyperbolicity of $\mathbb B(a,N).$
The desired inclusion now follows.
Next, we let $\Om'_r$ be the connected component of $\Om_r$ that contains 
$U_{\Om'_N} (z_{N'},4\eta_{N'}).$
By Lemma 3.6 (b), there exists a constant $C'>1$ such that for  $m \ge N', l \ge N'$ we have
$$k_{\Om'_r} (z_m, z_l) \le k_{U_{\Om'_N} (z_{N'}, 4\eta_{N'})}(z_m,z_l)\le C'k_{\Om'_N}(z_m,z_l).$$
This implies that $\{z_n\}_{n \ge N'} \subset \Om'_r$ is a $k_{\Om'_r}-$Cauchy sequence, giving a contradiction.
\vk.3cm 
\no $(d)\Rightarrow (a).$ Let $\{z_n\}$ be a $k_{\Omega}$-Cauchy sequence of $\Omega$.
We must show that $\{z_n\}$ is convergent to some point in $\Om$ in the original topology of $E$.
For $n \ge 1$, we set 
$$\rho'_n:=\sup\limits_{m, l\geq n}k_{\Omega}(z_m,z_l)<+\infty.$$
Then $\rho'_n \downarrow 0$. We also let
$$U_\Om (z_n, 4\rho'_n):=\{z\in \Omega: k_{\Omega}(z_n, z)<4\rho'_n\}.$$
We first claim that there exists $N''$ such that $U_\Om (z_{N''}, 4\rho'_{N''})$ is bounded.
If the claim is false, then there exists a sequence $\{w_n\} \subset \Om$ 
such that
$$\Vert w_n\Vert>n, k_\Om (z_n,w_n)<4\rho'_n.$$
Using the triangle inequality we deduce that $w_n \to \infty$ is also a $k_\Om-$Cauchy sequence. This violates the assumption that $\infty$ is a $k'-$point of $\Om.$ The claim follows.
Now, we 
apply Lemma 3.6 (b) to find a constant $C''>1$ such that
$$k_{U_\Om (z_{N''}, 4\rho'_{N''})}(p, q)\leq C''k_{\Omega}(p, q) \ \forall p, q\in U_\Om (z_{N''},\rho'_{N''}).$$
Hence $\{z_n\}_{n \ge N''}$ is a $k_{U_\Om (z_{N''}, 4\rho'_{N''})}$-Cauchy sequence. 
Choose
$R>0$ so large such that 
$$U_\Om (z_{N''},4\rho'_{N''}) \subset \mathbb B(0,R).$$
Observe that $\mathbb B(0,R)$ is complete hyperbolic, so
$z_n\to p\in \overline {\Omega}$ (in the original topology of $E$). Since every (finite) boundary point of $\Om$ is a 
$k'$-point, we have, in fact, $p\in \Omega$. Hence $\Omega$ is $k_{\Omega}$-complete.
\vskip0,3cm
\noindent
The above theorem yields the following partial generalization of Theorem 1 in \cite{G}.
\vk0,2cm
\noindent
{\bf Corollary 3.7.} {\it Let $\Om$ be a unbounded domain in a Banach space $E$. Then $\Om$ is complete hyperbolic if 
the following conditions are satisfied:

\noindent
(a) Every $a \in \partial \Om$ admits a local peak holomorphic functions i.e., there exists a \nhd\ $U$ of $a$ and a holomorphic function $h_a$ such that $\vert h_a\vert<1$ on $U \cap \Om$ and $\lim_{z \to a} \vert h_a (z)\vert=1.$

\noindent
(b) There exists a peak holomorphic function at $\infty$ for $\Om$ i.e, there exists a holomorphic function $h$ on 
$\Om$ such that $\vert h\vert<1$ on $\Om$ and 
$\vert h(z)\vert \to 1$ as $z \to \infty.$}
\vskip0,2cm
\noindent
{\it Proof.} We first show that $\infty$ is a $k'-$point for $\Om.$
Indeed, by the assumption we see that $h$ maps $\Om$ into $\De.$ 
Moreover, for each sequence $\{a_n\} \subset \Om, a_n \to \infty$ we may select a subsequence $\{a_{n_k}\}$ such that
$h(a_{n_k}) \to a' \in \partial \De.$ Since $a'$ is a $k'-$point of $\De$, by Proposition 2.6 we infer that $a$ is a $k'-$ point of $\Om.$ By the same reasoning, every finite boundary point of $\Om$ is also a $k'-$point.
Thus, we may apply Theorem 3.5 to complete the proof.
\vskip0,2cm
\noindent
{\bf Remark.} Let $\mathbb B^\infty$ be the unit ball in $l^\infty$ and
$\{a_n\}_{n \ge 1}$ be a sequence of positive numbers such that $\sum_{n \ge 1} a_n<\infty.$
Consider the open set
$$\Om:=\{z=(z_0, z_1, \cdots, z_n,\cdots) \in \mathbb C \times \mathbb B^\infty \subset l^\infty: \rho(z):=\Im z_0+\sum_{n \ge 1} a_n\vert z_n\vert^2<0\}.$$
It is easy to check that $\rho$ is a convex function on $l^\infty$. This implies that 
$\Om$ is convex. Thus, every $a \in \partial \Om$ admits a local peak holomorphic function.
Let $$f(z_0,z')=\frac{z_0+i}{z_0-i}, \  (z_0, z') \in l^\infty, z_0 \ne i.$$
We then have
$$
\begin{aligned}
\vert f(z_0, z')\vert^2-1&=\big \vert \frac{z_0+i}{z_0-i}\big \vert^2-1\\
&=\frac{4\Im z_0}{\vert z_0-i\vert^2} \le \frac{4\Im z_0}{\vert \Im z_0-1\vert^2}, \ \forall z_0 \ne i.
\end{aligned}$$
It then follows, since $\Im z_0<0$ on $\Om,$ that $\vert f\vert<1$ on $\Om$.
Moreover, by some easy estimates we also obtain $\vert f(z)\vert \to 1$ as $z \to \infty, z \in \Om.$
Hence $f$ is a peak holomorphic function at $\infty$ for $\Om.$ By Theorem 3.5 we conclude that $\Om$ is complete hyperbolic.
\vskip0.3cm 
\no
The next result roughly says that SVP and local SVP are equivalent for domains having $\infty$ as a $t-$point.
\vskip0,2cm
\noindent
{\bf Theorem 3.8}. {\it Let $\Omega$ be a unbounded domain in $E$ such that $\infty$ is  a $t$-point of $\Om$. Then the following statements are equivalent:

\noindent
(a) $\Omega$ has SVP.

\noindent
(b) For every $p\in \partial \Omega$, there exists a neighborhood $U$ of $p$ such that $U\cap \Omega$ has SVP.}
\vk0,2cm
\noindent
Since SVP and tautness are the same in $\mathbb C^n$, Theorem 3.8 essentially generalizes Proposition 3.2 in \cite{NP} 
(see also Proposition 2 in \cite{G}).
\vk.3cm 
\no{\it Proof}. $(a) \Rightarrow (b)$ is trivial. 

\noindent $(b) \Rightarrow (a).$
Let $A$ be a connected open subset of a Banach space and $\{f_j\}_{j \ge 1} \subset \Hol(A, \Omega)$ be a sequence such that $Z^u_{\{f_j\}}\ne\emptyset$. 
We first claim that $\{f_j\}$ is locally uniformly bounded on $A$. Indeed, if this is not so then 
$\exists A \ni z_j\to z_0\in A$ such that $f_j(z_j)\to \infty$.
Fix $\xi \in Z_{\{f_j\}}.$ 
Then, by the decreasing property of the Lempert functions we obtain
$$l_\Om (f_j(z_j), f_j(\xi)) \le l_{A} (z_j, \xi).$$
Since $\infty$ is a $t-$ point of $\Om,$ by letting $j \to \infty$ and taking limsup, we obtain
$$1=\limsup\limits_{j\to\infty} l_{A} (z_j, \xi) \le l_{A} (z_0, \xi) <1,$$
which is absurd.
Here the middle inequality follows from upper semicontinuity of $l_A$ (cf. Lemma 2.2). 

Thus we have shown that $\{f_j\}$ is locally uniformly bounded on $A$.
Now, we use Lemma 3.3 to get that $f_j\rightarrow f$ in $\Hol(A,E)$. Put $A':=f^{-1}(\Omega)$. Then, obviously  
$A' \ne \emptyset$ and $A'$ is open. Suppose that $A' \ne A$. Then, we can find 
$z_0\in A \cap \partial A'$. By (ii), we may choose a neighborhood $V$ of $f(z_0)$ such that $V\cap \Omega$ has SVP. Since $f_j\to f$ in $\Hol(A,E)$ there exists a neighborhood $U$ of $z_0$ and $j_0\geq 1$ such that 
$$f_j(U)\subset V\cap \Omega \ \forall j\geq j_0.$$
It follows that $f_j\to f$ in $\Hol(A, V\cap\Omega)$, because $V\cap \Omega$ has SVP and 
$A'\cap U\subset Z^u_{\{f_j|_U\}}$. Hence $A'= A,$ a contradiction.
Summing up, we have proved our claim that $\{f_j\}$ is convergent in $\Hol(A,\Omega)$.
This is the desired conclusion.
\vk.3cm 
\no
As an illustration of the above theorem we have the following equivalence between Vitali properties on certain classes
of unbounded domains in Banach spaces.
\vk0,2cm
\noindent
{\bf Corollary 3.9.} {\it Let $\Omega$ be a unbounded domain in a Banach space such that $\infty$ is a $t$-point. Then 
$\Omega$ has SVP if only if $\Omega$ has WVP.}
\vk.3cm 
\no{\it Proof.}  Assume $\Om$ has WVP. Fix $a \in \partial \Om$ and a ball $\mathbb B$ around $a$. We claim that
$\Om \cap \mathbb B$ has SVP. For this, we let $A$ be a connected open subset of a Banach space and $\{f_j\}$ be  sequence in $\Hol(A, \Om \cap \mathbb B)$ such that
$Z^u_{\{f_j\}} \ne \emptyset.$ 
By Lemma 3.3, we infer that $\{f_j\}$ is convergent to $f \in \Hol(A,E)$. Set $U:=f^{-1} (\Om \cap \mathbb B).$
Then $U$ is open and non-empty since $Z_{\{f_j\}} \subset U.$ Notice that $U \subset Z_{\{f_j\}} \cap Z^u_{\{f_j\}}.$ 
Hence, as $\Om$ has WVP, $\{f_j\}$ is convergent to $f \in \Hol(A,\Om).$
It remains to check that $f(A) \cap \mathbb B =\emptyset.$ Assume otherwise, then there exists 
$\xi \in A$ such that $f_j (\xi) \to \partial \mathbb B.$
Fix $a \in Z_{\{f_j\}}.$
By the decreasing  property of Kobayashi pseudo-distance we obtain
$$\sup_{j \ge 1} k_{\mathbb B} (f_j(\xi), f_j(a)) \le k_A (\xi, a)<\infty.$$
We arrive at a contradiction, since $\mathbb B$ is complete hyperbolic.

\noindent
{\bf Remark.} In light of the above result, the following question is of interest to us: Let $\Om$ be an unbounded  domain in a Banach space having SVP. Is $\infty$ necessarily a $t-$point of $\Om$?
\vk0.2cm
\noindent
The result below provide explicit examples of unbounded domains in Banach spaces having SVP.
\vk.3cm 
\no{\bf Corollary 3.10.} {\it Let $\Omega$ be a unbounded domain in a Banach space $E$ such that $\infty$ is a $t$-point. Then $\Omega$ has SVP if one of the following holds:

\no (i) Every (finite) boundary point of $\Omega$ is a $t$-point;

\no (ii) Every (finite) boundary point of $\Omega$ admits a barrier.}
\vk.3cm 
\no{\it Proof}. Assume that (i) holds. 
Fix $p\in \partial \Omega$ and $r>0$. We claim that
$U: =\mathbb B(p,r) \cap \Om$ has SVP. Indeed, let $A$ be a connected open subset of a Banach space and  
$\{f_j\} \subset\Hol(A,U)$ with $Z^u_{\{f_j\}}\ne\emptyset$. 
Since $\mathbb B(p, r)$ is bounded in $E$, by Lemma 3.3, 
$f_j\to f$ in $\Hol(A,E)$. 
Now we suppose that there exists $\lambda \in A$ such that 
$$q:=f(\lambda)=\lim_{j \to \infty} f_j(\lambda)\in \partial U.$$
Since, by Proposition 2.4 every boundary point of $\mathbb B(p,r)$ is a $t-$point of $U$ and
since every boundary point of $\Om$ is a $t-$point of $\Om$ by the assumption, we conclude that $q$ is a $t-$point of $U.$
Next, we choose $\la' \in A \cap Z_{\{f_j\}}$. Then, by the decreasing property of the Lempert functions we obtain
$$1=\sup_{j \ge 1} l_U (f_j(\la), f_j (\la')) \le l_A (\la, \la')<1,$$
which is absurd. It follows that $f(A)\subset U$.
Hence $f_j\to f$ in $\Hol(A,U)$.
This is exactly our claim. It remains to apply Theorem 3.8 to conclude that $\Om$ has SVP.

Finally, suppose that $(ii)$ is true. Then, we first apply Theorem 3.1 to see that $\Omega$ is hyperbolic. Next, we use Proposition 2.4 to get that every finite boundary point of $\Om$ is a $t$-point. By (i) $\Omega$ has SVP. This fully complete our proof.

\noindent
{\bf Remark.} Let $E$ be a Banach space, $\mathbb B$ be the unit ball in $E$
and $u \in PSH(\mathbb B) \cap \mathcal C (\bar {\mathbb B})$ be such that $u \ge 0$ on $\mathbb B$.
Set $E':=\mathbb C \times E.$
Consider the unbounded domain
$$\Om:=\{z=(z_0,z')\in E': z' \in \mathbb B, v(z):=\Im z_0+u(z')<0\}.$$
Since $\Om$ is contained in the product $\{\Im z_0<0\} \times \mathbb B$ of hyperbolic domains in $\mathbb C$ and in
$E$ respectively, we infer that $\Om$ is hyperbolic as well.
We now claim that every $a \in \partial \Om \cup \{\infty\}$ admits a barrier. Indeed, for $a \in \partial \Om,$ we may pick $\max \{v(z), \Vert z'\Vert^2-1\}$ as a barrier at $a.$ For $a=\infty$, by the same reasoning as in the remark that follows we see that 
$$\va (z):=\big \vert \frac{z_0+i}{z_0-i}\big \vert^2-1$$
is a \psh\ barrier at $\infty.$
This proves our claim. Hence we may apply Corollary 3.10 (b) to conclude that $\Om$ is an unbounded domain having SVP.
\vk.3cm
\no{\bf Proposition 3.11.} {\it Let $\theta: X \to Y$ be a holomorphic map between connected open subsets of Banach spaces. Assume that $Y$ has SVP (resp. WVP) and for every $y\in Y$ there exists a connected hyperbolic \nhd\ $U$ of $y$ such that $\theta^{-1}(U)$ has SVP (resp. WVP). Then $X$ has  SVP (resp. WVP).}
\vk.3cm 
\no{\it Proof}. For the ease of exposition we only deal with the case where $Y$ and all the fibers $\theta^{-1} (U)$ have
SVP. The other case can be treated analogously.
Fix a connected open subset $A$ in a Banach space and $\{f_j\}\subset \Hol(A,X)$ with $Z^u_{\{f_j\}}\ne\emptyset$. 
Put $g_j:=\theta\circ f_j$. Since $Z^u_{\{g_j\}}\supset Z^u_{\{f_j\}}\ne\emptyset$ and $Y$ has SVP we have $g_j\to g$ in $\Hol(A, Y)$. Take $\lambda_0\in Z^u_{\{f_j\}}$ and $U\ni g(\lambda_0)$ such that $\theta^{-1}(U)$ has SVP. Since $U$ is hyperbolic and since $g_j\to g$  in $\Hol(A,Y)$, there exists a neighborhood $V\ni \lambda_0$ and $j_0 \ge 1$ such that
$$g_j(V)\subset U\ \forall j \ge j_0.$$
It follows that 
$$f_j(V)\subset \theta^{-1}(U)\ \forall j \ge j_0.$$
Since $Z^u_{\{f_j|_V\}}\ni \lambda_0$ we have $\{f_j|_V\}$ converges in $\Hol(U,X)$. Put
$$\Omega=\bigcup\Big\{V \subset A: V \ \text{is open}\ \{f_j|_V\} \ \text{is convergent in}\ \Hol(V,X)  \Big\}.$$ 
Then $\Omega$ is open and nonempty. Let $\lambda_1\in\partial \Omega$. The above reasoning implies that there exist open \nhd s 
$U' \ni g(\lambda_1)$ and $V' \ni \lambda_1$ such that
$$f_j(V')\subset U' \ \forall j \ge j_0$$
and $\theta^{-1}(U')$ has SVP.
Hence $\{f_j|_{V'}\}$ converges in $\Hol(V',X)$. This implies $\lambda_1\in\Omega$. Consequently $\Omega=A$.
This finishes our proof.
\vk.5cm 
\no{\bf Remarks.} (a) By the same proof as in \cite{Ko}, we can prove an analogous result to Proposition 3.10 where Vitali properties is replace by complete hyperbolicity.

\noindent
(b) In the case where $X, Y$ are complex manifolds, an analogous result to Proposition 3.11, where Vitali properties is replaced by taut property, is obtained in Proposition 5.1.8 of \cite{Ko}. Notice that, it was assumed there that the holomorphic map $\theta$ is {\it proper}. 
Since Vitali properties and tautness coincide in the case of complex manifolds 
(see Theorem 3.14 in \cite{KDK}), Proposition 3.11 is somewhat stronger than the above mentioned result in \cite{Ko}.

\noindent
(c) Let $\Om$ be a a domain having $SVP$ in a Banach space and $\va$ be a continuous \psh\ function on $\Om.$
Consider the Hartogs domain $\Om' \subset \Om \times \mathbb C$ defined by
$$\Om':=\{(z,w): z \in \Om, \log \vert w\vert+\va (z)<0\}.$$
Then we have the projection map $\theta: \Om' \to \Om, (z,w) \mapsto z.$
Fix $z_0 \in \Om,$ we let $U \subset \Om$ be a small ball of radius $r>0$ around $z_0$. It follows that
$$U':=\theta^{-1} (U)=\{(z,w) \in U \times \mathbb C: \psi (z,w):=\log \vert w\vert+\va(z)<0.\}$$
Then each boundary point $\partial U'$ admits the barrier $\max\{\psi (z,w), \Vert z-z_0\Vert-r\}.$
Then, by the {\it proof} of Corollary 3.10 we conclude that $U'$ has SVP. Thus, $\Om'$ has SVP by Proposition 3.11.
\vk.5cm

\vk.3cm

\no Department of Mathematics, Hanoi National University of Education. 

\no Address: 136 Xuan Thuy Street, Hanoi, Vietnam.

\no Email-address: dieu\_vn@yahoo.com (Nguyen Quang Dieu).

\no Email-address: nvkhiem@hnue.edu.vn (Nguyen Van Khiem).

\no College of Education, Trung Trac, Vinh Phuc.

\no Email-address: thanhhungcdsp@gmail.com (Le Thanh Hung).

\end{document}